\documentclass[12pt]{article}
\usepackage{amsmath}
\usepackage{graphicx}
\usepackage{amscd}
\usepackage{amsfonts}
\usepackage{amssymb}
\newtheorem{theorem}{Theorem}[section]
\newtheorem{proposition}[theorem]{Proposition}

\newtheorem{corollary}[theorem]{Corollary}

\begin{document}
\title{Modular deformations of analytic polyhedra\footnote{This research was
supported by a grant from the GIF, German-Israeli Foundation of
Scientific Research and Development.}}
\author{V. P. Palamodov}
\date{20.06.2005}
\maketitle

\section{Introduction}

The notion of modular deformation helps to define a canonical
analytic structure in the category of analytic objects. In this
paper, we discuss modular deformations of complex analytic spaces.
The category of germs of complex spaces is not appropriate for
this purpose since the fiber of a deformation of a germ is no more
a germ. Instead, we consider the category of analytic polyhedrons
(see Sec.4). Let $f:\mathcal{X}\rightarrow S$ be a deformation of
an analytic space or of a polyhedron. Take a point $\circ$ in the
base $S$ and consider the pair $\left(  S,\circ\right)  $ as a
germ of complex space. A \textit{modular} stratum is an analytic
subspace $\left( M,\circ\right)  \subset\left(  S,\circ\right) $
such that the uniqueness property holds: any germ morphisms
$g:R\rightarrow\left(  S,\circ\right)  $ and $h:R\rightarrow\left(
M,\circ\right)  $ induce isomorphic deformations
$f\times_{S}g\cong f\times_{S}h$ only if $g=h.$ The
\textit{modular} deformation is the restriction of $f$ to $M,$
that is the deformation $f_{M}\doteq f\times_{S}M.$ Suppose that
$f:\mathcal{X}\rightarrow\left(
S,\circ\right)  $ is a versal deformation of the fiber $X_{\circ}%
=f^{-1}\left(  \circ\right)  $ and $M$ is the \textit{maximal}
modular stratum. It has the special feature: for any other versal
deformation $f^{\prime}:\mathcal{X}^{\prime}\rightarrow
S^{\prime}$ of $X$ and an isomorphism $i:S\rightarrow S^{\prime}$
such that $i\times f^{\prime}\cong f$, the restriction of $i$ to
$M$ is uniquely defined and $M^{\prime}\doteq i\left(  M\right)  $
is the maximal modular stratum in $S^{\prime}$. The isomorphism
$\left(  i\times f^{\prime}\right)  _{M}\cong f_{M}$ is also
canonically defined. Moreover, the maximal modular stratum $M$ (if
it exists) possesses the semi-local property: for any point $s\in
M$ close to the marked point $\circ$ the germ $\left(  M,s\right)
$ is also maximal modular for the fiber $X_{s}\doteq f^{-1}\left(
s\right)  $. Therefore two versal deformations $f$ and
$f^{\prime}$ can be glued together (amalgamated) along its modular
strata, say $M$ and $M^{\prime},$ provided there is given an
isomorphism $\phi:f^{-1}\left(  s\right)
\tilde{\rightarrow}f^{\prime -1}\left(  s^{\prime}\right)  $ for
$s\in M,s^{\prime}\in M^{\prime}.$ This gives rise to a flat
morphism $F:X\rightarrow M\sqcup M^{\prime}$ where $M\sqcup
M^{\prime}$ denotes the amalgam over some open neighborhood of the
points $s$ and $s^{\prime}.$ Repeating this construction, we can
get a flat morphism of complex spaces
$F:\mathcal{X}\rightarrow\mathcal{M}$ which is a maximal modular
deformation of each its fiber. This conception is close to that of
fine moduli space in the sense of Mumford, \cite{Mum}. For each
point $s\in\mathcal{M},$ the automorphism group of the space
$X_{s} $ acts on $\mathcal{M}$ which may possible further
amalgamation and the global amalgam may have complicated
structure.

We formulate here some general results on existence and properties
of maximal modular deformations in the category of complex
analytic polyhedra and discuss several examples of modular
deformation of isolated singularities. V.Arnold's classification
table of hypersurface singularities \cite{AGV} is a rich source of
such examples. For many examples of complete intersection
singularities modular deformations are contained in the paper of
A. Aleksandrov \cite{Al1}, further study see in \cite{Al}.

Several complicated examples, which have been beyond the reach of
`by hand' calculations, were computed and studied by B.Martin
\cite{Ma2} and by T. Hirsch and B. Martin \cite{HM}. These authors
applied a specialized computer algebra program based on SINGULAR.
We shall see how these new local modular families glue together
and global modular families emerge. In several cases the base of a
modular family can be made compact by gluing together sufficiently
many local modular deformations.

I thank Bernd Martin for many helpful discussions.

\section{Rudiments of the deformation theory}

Remind some basic definitions. A germ of complex analytic spaces (which will
be called simply `germ') is a pair $\left(  Z,\mathcal{O}(Z)\right)  $, where
$Z$ is the germ of complex analytic set in $\left(  \mathbb{C}^{n}%
,\circ\right)  $ for some $n$ and $\mathcal{O}(Z)$ is the coherent sheaf of
analytic $\mathbb{C}$-algebras on $\left(  \mathbb{C}^{n},\circ\right)  $ such
that $\mathrm{supp}\,\mathcal{O}\left(  Z\right)  =Z$. We use the notation
$\circ$ for the marked point in $\mathbb{C}^{n}$ and in $Z$. Any morphism of
germs $\left(  W,\mathcal{O}(W)\right)  \rightarrow\left(  Z,\mathcal{O}%
(Z)\right)  $ is a pair $\left(  f,\phi\right)  $ where $f:\left(
W,\circ\right)  \rightarrow\left(  Z,\circ\right)  $ is a mapping of germs of
analytic sets and $\phi:f^{\ast}\left(  \mathcal{O}(Z)\right)  \rightarrow
\mathcal{O}(W)$ is a morphism of sheaves of analytic $\mathbb{C}$-algebras
over $W$. Restricting $\phi$ to the marked point yields the morphism of
analytic $\mathbb{C}$-algebras $\phi_{\circ}:\,\mathcal{O}_{\circ
}(Z)\rightarrow\mathcal{O}_{\circ}(W)$. The morphism $(f,\phi)$ is called
embedding, if $\phi_{\circ}$ is a surjection. Vice versa, for any analytic
$\mathbb{C}$-algebra $A$ there exists and is uniquely defined the germ
$\left(  Z,\mathcal{O}\left(  Z\right)  \right)  $ such that $A\cong
\mathcal{O}_{\circ}\left(  Z\right)  $. In particular, the algebra
corresponding to the simple point $\circ$ is the field $\mathbb{C}$. For an
arbitrary germ $Z$ there is the canonical embedding $\circ\rightarrow Z$; the
corresponding algebra homomorphism is the canonical morphism $\mathcal{O}%
_{\circ}(Z)\rightarrow\mathbb{C}$. The kernel of this morphism is the maximal
ideal in the algebra $A=\mathcal{O}_{\circ}(Z);$ it is denoted $\mathfrak
{m}(Z)$. The germ $Z=\circ$ is called \textit{simple} point; a germ $Z$ is
called \textit{fat} point, if the corresponding algebra has a finite dimension
over $\mathbb{C}$.

The fiber product operation is well defined in the category of analytic
spaces; it is denoted $f\times_{X}g$ or $Y\times_{X}Z$ for germ morphism
$f:Y\rightarrow X,\;g:Z\rightarrow X.$ If $g$ is an embedding, we call the
product $f\times_{X}g$ restriction $f$ to $Z$ and denote it $\left.  f\right|
Z$. Let $V$ be a $\mathbb{C}$-vector space of finite dimension; we denote by
$\mathbb{D}\left[  V\right]  $ the $\mathbb{C}$-algebra that is isomorphic to
$\mathbb{C}\oplus V$ as vector space with the multiplication rule
$(a,u)\cdot(b,v)=(ab,av+bu)$. This is an analytic algebra of a fat point
denoted $P\left(  V\right)  $. All germs of this kind form a category denoted
$\emph{D}$, where for any objects $P\left(  U\right)  ,P\left(  V\right)  $
the set $\mathrm{Hom}(P\left(  U\right)  ,P\left(  V\right)  )$ is bijective
to $\mathrm{Hom}_{\mathbb{C}}(U,V).$ This induces a natural structure of
vector space over $\mathbb{C}$ in the set $\mathrm{Hom}(P\left(  U\right)
,P\left(  V\right)  )$; the composition of morphisms is a bilinear operation.
For any germ $Z$ there is defined the germ $T(Z)\in\emph{D}$ and the embedding
$t(Z):T(Z)\rightarrow Z$ such that $\mathcal{O}(T(Z))=\mathcal{O}%
(Z)/\mathfrak{m}^{2}(Z).$ The germ $T(Z)$ is called the \textit{tangent space}
of the germ $Z$.\newline \textbf{Definition 1. }Let $X$ be a complex analytic
space. A \textit{deformation } of $X$ over a germ $\left(  Z,\circ\right)  $
is a pair $(f,i)$, where $f:\mathcal{X}\rightarrow\left(  Z,\circ\right)  $ is
a flat morphism and $i:X\rightarrow f^{-1}\left(  \circ\right)  =f\times\circ$
is an isomorphism of complex analytic spaces. For any morphism of germs
$h:W\rightarrow Z$ and a deformation $\left(  f,i\right)  $ of $X$ over $Z$
the fiber product $f\times_{Z}h:\mathcal{X}\times_{Z}W\rightarrow W$ is also
flat. The pair $\left(  f\times_{Z}h,ji\right)  $ is a deformation of $X$ with
base $W,$ where $j:f^{-1}\left(  \circ\right)  \cong\left(  f\times
_{Z}h\right)  ^{-1}\left(  \circ\right)  $ is the canonical bijection. In the
same way, one can treat deformations of $\mathbb{Z}_{2}$-graded, $\mathbb{Z}%
$-graded analytic spaces, deformations of germs, of fiber bundles, of coherent
sheaves and the like.

A \textit{versal }deformation of the space $X$ is a deformation $(f,i),$
$f:\mathcal{X}\rightarrow S$ such that for any deformation $\left(
g,j\right)  ,g:\mathcal{Y}\rightarrow R$ there exists a morphism of germs
$h:R\rightarrow S$ and a isomorphism $\gamma:\mathcal{X}\times_{S}%
R\rightarrow\mathcal{Y}$ over $R$ such that the diagram commutes
\[%
\begin{tabular}
[c]{lllll}%
$\mathcal{X}\times_{S}\circ$ &  & $\overset{\gamma\times_{R}\circ}%
{\rightarrow}$ &  & $\mathcal{Y}\times_{R}\circ$\\
& $\overset{i}{\nwarrow}$ &  & $\overset{j}{\nearrow}$ & \\
&  & $X$ &  &
\end{tabular}
\]
A pair $(f,i)$ is called \textit{universal} deformation, if the morphism $h$
is unique. The above definition can be applied for the category $\emph{D}$.
Suppose that a space $X$ possesses a universal deformation $\delta
:\mathcal{X}\rightarrow T_{X}$ in the category $\emph{D}$. Take an arbitrary
germ $S$ and a deformation $\left(  f,i\right)  $ of $X$ over $S$; consider
the morphism $f_{T}=$ $f\times_{S}T_{\circ}\left(  S\right)  $. It is a
deformation of $X$ with the base $T_{\circ}\left(  S\right)  \in\emph{D}$. Due
to the universality property of $\delta$ a morphism $\mathrm{D}_{\circ
}f:T_{\circ}(S)\rightarrow T_{X}\;$is defined in the category $\emph{D}$ such
that $\mathrm{D}_{\circ}f\times\delta\cong f_{T}.$ It is called
\textit{Kodaira-Spencer} mapping; this mapping is linear since it belongs to
$\emph{D}$. Let $\left(  f,i\right)  $ be a versal deformation of $X$ with a
base $S.$ There exists a morphism $t:T_{X}\rightarrow S$ such that
$f\times_{S}t\cong\delta.$ Then $\mathrm{D}_{\circ}f\cdot t=\operatorname{id}$
since $\delta$ is universal, consequently $\mathrm{D}_{\circ}f$ is surjective.
The pair $\left(  f,i\right)  $ is called \textit{miniversal }(or minimal
versal), if $\mathrm{D}_{\circ}f$ is a injection, hence a bijection.\newline
\textbf{Definition 2.} \cite{P2}, \cite{P3} Let $f$ be a deformation of $X$ as
above. A subgerm (stratum) $M\subset S$ is called \textit{modular}, if for any
morphisms of germs $h:R\rightarrow M,\,g:R\rightarrow S$ the equation
$f\times_{S}h=f\times_{S}g$ implies that $h=g$. A modular stratum $M\subset S$
is called \textit{maximal}, if it contains any other modular stratum. If
$f:\mathcal{X}\rightarrow S$ is a miniversal deformation and $M$ is the
maximal modular stratum, the restriction $f_{M}$ is called maximal modular
deformation of $X.$ Any universal deformation is, of course, maximal modular.

A similar conception (prorepresenting stratum) was introduced in \cite{La} for
the formal deformation theory of affine schemes. See \cite{LP} for further
results. The notion of modular deformation was also treated in \cite{KSt} in a
more general setting.

\begin{proposition}
\label{p2}Let $f:\mathcal{X}\rightarrow S$ and $f^{\prime}:\mathcal{X}%
^{\prime}\rightarrow S^{\prime}$ be miniversal deformations of a space $X$ and
$M,\,M^{\prime}$ are respective maximal modular strata. There exists a
isomorphism $h:M^{\prime}\rightarrow M$ such that $\left.  f^{\prime}\right|
M^{\prime}\cong\left.  f\right|  M\times h.$ The morphism $h$ is unique.
\end{proposition}

$\blacktriangleleft$ There exist morphisms $g:S^{\prime}\rightarrow S$ and
$g^{\prime}:S\rightarrow S^{\prime}\;$such that $f^{\prime}\times g^{\prime
}\cong f,\;f\times g\cong f^{\prime}.$ We have $\mathrm{D}_{\circ}%
f\cdot\mathrm{d}g=\mathrm{D}_{\circ}f^{\prime}$ and $\mathrm{D}_{\circ
}f^{\prime}\cdot\,\mathrm{d}g^{\prime}=\mathrm{D}_{\circ}f.$ The differentials
$\mathrm{d}g,$ $\mathrm{d}g^{\prime}$ are inverse one to another, since
$\mathrm{D}_{\circ}f,\;\mathrm{D}_{\circ}f^{\prime}$ are bijections. Therefore
$g$ is a isomorphism of germs and the mapping $g:M^{\prime}\rightarrow S$ is a
modular germ. It is factorized through a uniquely defined morphism
$h:M^{\prime}\rightarrow M,$ since $M$ is maximal modular.
$\blacktriangleright$

\section{Tangent cohomology and criterion of a modular stratum}

Let $g:X$ be a complex analytic space; the graded tangent sheaf $\oplus
_{q=0}\mathcal{T}^{q}(X)$ is defined on $X.$ It is a $\mathbb{Z}_{+}$-graded
sheaf algebra Lie, which means that there is defined a bracket operation that
satisfies the graded commutation and Jacobi identities. Moreover, this has a
natural structure of $\mathcal{O}\left(  X\right)  $-module which agrees with
the Lie algebra structure in a natural way. For any $q\geq0$ the sheaf
$\mathcal{T}^{q}$ is a coherent $\mathcal{O}\left(  X\right)  $-sheaf. In
particular, the term $\mathcal{T}^{0}\left(  X\right)  $ is the sheaf of
tangent fields on $X$.

The global tangent cohomology $T^{\ast}(X)=\sum_{0}^{\infty}T^{n}(X)$ is a
$\mathbb{Z}_{+}$-graded algebra and there is a spectral sequence $E^{\ast}$
that converges to the tangent cohomology $T^{\ast}(X)$ with the second term
$E_{2}^{pq}=H^{p}(X,\mathcal{T}^{q}(X))$, see \cite{P1}. The term
$T^{0}\left(  X\right)  $ is the Lie algebra of tangent fields on $X.$ If the
term $T^{1}\left(  X\right)  $ is of finite dimension, it represents the base
$T_{X}$ of the universal deformation $\delta$ of $X$ in the category
$\emph{D}$ as in Sec.2. If $g:X\rightarrow Y$ is a morphism of complex spaces,
a vertical tangent field $t$ on $g$ is a tangent field on $X$ with the
property $t\left(  g^{\ast}\left(  a\right)  \right)  =0$ for an arbitrary
$a\in\mathcal{O}\left(  Y\right)  .$ The notation $T^{0}\left(  X/Y\right)  $
means the space of vertical tangent fields. For any point $y\in Y$ the
restriction mapping $T^{0}\left(  X/Y\right)  \rightarrow T^{0}\left(
X_{y}\right)  $ is canonically defined where $X_{y}=g^{-1}\left(  y\right)  $
is the complex subspace of $X.$ There is a general criterion of modularity:

\begin{theorem}
\label{t1}Let $f:\mathcal{X}\rightarrow\left(  S,\circ\right)  $ be a
deformation of a complex analytic space $X$. Then \newline (i) the simple
point $\circ\in S$ is modular, if and only if the Kodaira-Spencer mapping
$\mathrm{D}_{\circ}f$ is injective;\newline (ii) if $\mathrm{D}_{\circ}f$ is
injective, then a subgerm $M\subset S$ is modular, if for any fat point
$Z\subset M$ the restriction mapping $T^{0}(\mathcal{X}\times_{M}%
Z/Z)\rightarrow T^{0}(X)$ is surjective; if $\mathrm{D}_{\circ}f$ is
bijective, this condition is necessary as well;\newline (iii) $T_{\circ}(M)$
coincides with the space of tangent vectors $t\in T_{\circ}\left(  S\right)  $
that satisfy the equation $[\mathrm{D}_{\circ}f(t),v]=0$ for any $v\in
T^{0}(X)$.
\end{theorem}

A proof is similar to that of \cite{P3}, Proposition 5.1.

\section{Analytic polyhedrons}

The existence of maximal modular deformation and some its properties were
stated in \cite{P3} for an arbitrary compact complex analytic space. Here we
formulate similar results for analytic polyhedra. First we recall some
definitions of \cite{P6}.\newline \textbf{Definition 3.} For an arbitrary
integer $n$ we fix a coordinate space $\mathbb{C}^{n}$, i.e. a complex vector
space with a marked system of linear coordinates $w_{1},\ldots,w_{n}$. Denote
by $D^{n}$ the closed unit polydisk in $\mathbb{C}^{n}$. A complex analytic
$n$-\textit{polyhedron }is a pair $(X,\varphi)$, where $X$ is a complex
analytic subspace of a complex space $\tilde{X}$ and $\varphi:\tilde
{X}\rightarrow\mathbb{C}^{n}$ is a holomorphic mapping (called a
\textit{barrier} map ) such that the set $\varphi^{-1}\left(  \bar{D}%
^{n}\right)  $ is compact and $X=\varphi^{-1}\left(  D^{n}\right)  .$ The
neighborhood $\tilde{X}$ of $X$ can be contracted, i.e. $X$ is thought as the
germ of a complex analytic space $\tilde{X}$ on the compact set $\bar
{X}=\varphi^{-1}\left(  \bar{D}^{n}\right)  .$ The set $\partial
X\doteq\varphi^{-1}\left(  \partial D^{n}\right)  $ is the boundary of the
polyhedron $X.$ A morphism of polyhedra $\left(  X,\varphi\right)
\rightarrow$ $\left(  Y,\psi\right)  $ is a pair $\left(  f,p\right)  ,$ where
$f:\bar{X}\rightarrow\bar{Y}$ of germs of complex spaces and $p:\mathbb{C}%
^{n}\rightarrow\mathbb{C}^{m}$ is a coordinate projection that make the
commutative diagram: $p\varphi=\psi f$.

Let $S$ be a complex analytic spaces; we call a \textit{relative }analytic
$n$-polyhedron (r.p.) over $S$ any pair $(X,\varphi)$, where $X$ is an open
subspace a complex space $\tilde{X}$ and $\varphi:\tilde{X}\rightarrow
\mathbb{C}^{n}\times S$ is a holomorphic mapping such that the set
$\varphi^{-1}\left(  \bar{D}^{n}\times S\right)  $ is proper over $S$ and
$\varphi^{-1}\left(  D^{n}\times S\right)  =X$. Let $(X,\varphi)$ be a r.p.
over $S$ and $h:R\rightarrow S$ be a morphism of complex spaces. The fiber
product $(X_{R},\varphi_{R})\doteq(X,\varphi)\times_{S}R$ is a r.p. over $R.$
In particular, for any point $s\in S$ the fiber of $(X,\varphi)$ is defined as
the product $(X_{s},\varphi_{s})=(X,\varphi)\times_{S}s$, which is an
(absolute) analytic polyhedron, called the fiber of $\left(  X,\varphi\right)  .$

A morphism $\left(  X,\varphi\right)  \rightarrow\left(  Y,\psi\right)  $ of
relative polyhedra over $S$ is any pair $\left(  f,p\right)  $ of holomorphic
mappings of germs $f:\bar{X}\rightarrow\bar{Y},\,p:\bar{D}^{n}\times
S\rightarrow\bar{D}^{m}\times S$ such that $p\varphi=\psi f$ and $p$ commutes
with the projections to $S.$ If $S=\circ$ this is a definition of morphism of
(absolute) polyhedra.\newline \textbf{Definition 4.} Let $S$ be a germ of
complex space with the marked point $\circ$, and $(X,\varphi)$ be a
polyhedron; we call a \textit{deformation }of this polyhedron with the base
$S$ any triple $(\mathcal{X},\Phi;\theta)$, where \newline (i) $(\mathcal{X}%
,\Phi)$ is a relative polyhedron over $S$, \newline (ii) $(\mathcal{X},\Phi)$
is flat over $S,$ i.e. the composition $\pi\Phi:\mathcal{X}\rightarrow
\mathbb{C}^{n}\times S\rightarrow S$ is a flat morphism of complex analytic
spaces, where $\pi$ is the projection. \newline (iii) $\theta:(X,\varphi
)\rightarrow(\mathcal{X},\Phi)\times_{S}\left\{  \circ\right\}  $ is an
isomorphism of polyhedra. A deformation $(\mathcal{X},\Phi;\theta)$ of a
polyhedron $(X,\varphi)$ is called \textit{versal}, if for any deformation
$(\mathcal{Y},\Psi;\eta)$ of the same polyhedron with a base $R$ there exist a
morphism of germs $h:R\rightarrow S$ and a isomorphism $\alpha:(\mathcal{Y}%
,\Psi)\rightarrow(\mathcal{X},\Phi)\times_{S}R$ such that $\left(
\alpha\times\circ\right)  \theta=\eta$. Definitions of modular stratum,
maximal modular deformation of an analytic polyhedron are given in the same
lines as in Definition 2.

\section{Versal and modular deformations}

Let $\mathbb{C}^{n}$ be again a coordinate space and $D^{n}$ be the closed
unit polydisk. For an arbitrary coordinate projection $p:\mathbb{C}%
^{n}\rightarrow\mathbb{C}^{m}$ we have $p(D^{n})=D^{m}$ and the sheaf
$\mathcal{O}(D^{n})$ is endowed with a structure of module over the sheaf
algebra $\mathcal{O}(D^{m})$.

\begin{theorem}
\label{2.1}Let $(X,\varphi)$ be an analytic $n$-polyhedron such that
\newline (i) for any coordinate projection $p:\mathbb{C}^{n}\rightarrow
\mathbb{C}^{m}$ the equations hold:
\begin{equation}
Tor_{k}^{\mathcal{O}(D^{m}}(\mathcal{T}_{\varphi}^{k+1},\mathbb{C)}%
=0,\;k=0,1,\dots,m\label{2.2}%
\end{equation}
in each point of the set $(\partial D)^{m}\times D^{n-m}$, where
$\mathcal{T}^{\ast}\doteq\mathcal{T}^{\ast}(X)$ and \newline (ii)
$\phi$ is finite over $\partial D^{n}.$ The polyhedron $\left(
X,\varphi\right)  $ has a miniversal deformation
$F:\mathcal{X}\rightarrow\left(  S,0\right)  $ where $S\cong
g^{-1}(0)$ and $g:\left(  T^{1}\left(  X\right)  ,0\right)
\rightarrow T^{2}(X)$ is a holomorphic mapping of germs such that
$g(0)=0,\,\mathrm{d}g(0)=0$.
\end{theorem}

\textbf{Remark}. The conditions (\ref{2.2}) are fulfilled for the trivial
barrier $\varphi=0$, since $\mathcal{O}(D^{0}\mathbb{)=C}.$ This implies the
existence of miniversal deformation for any compact complex space. This case
was studied by H.Grauert and other authors, see the survey \cite{P1}. Theorem
\ref{2.1} was proved in \cite{P6}. Note that the space $X$ that satisfies
(\ref{2.2}) need not to be smooth at the boundary of the polyhedron. The
condition (\ref{2.2}) for $k=0$ means that the sheaf $\mathcal{T}_{\varphi
}^{1}$ vanishes on the boundary of $D^{n}$, hence the support of this sheaf is
a compact subset of the interior of $D^{n}$. The spectral sequence $E_{2}%
^{pq}=H^{p}(X,\mathcal{T}^{q})\Rightarrow T^{\ast}\left(  X\right)  $ yields
\[
0\rightarrow H^{1}\left(  X,\mathcal{T}^{0}\right)  \rightarrow T^{1}%
(X)\rightarrow\Gamma(X,\mathcal{T}^{1})\rightarrow H^{2}\left(  X,\mathcal{T}%
^{0}\right)  \rightarrow...
\]
The term $\Gamma\left(  X,\mathcal{T}^{1}\right)  $ is of finite dimension
since $\mathcal{T}^{1}$ is supported by a compact set. From the Leray sequence
we find the exact sequence
\begin{equation}
H^{1}\left(  X,\mathcal{T}^{0}\right) \rightarrow H^{0}\left(  D^{n}%
,\mathcal{R}^{1}\varphi\left(  \mathcal{T}^{0}\right)  \right)  \overset
{d_{2}}{\rightarrow} H^{2}\left(  D^{n},\mathcal{R}^{0}\varphi\left(
\mathcal{T}^{0}\right)  \right)  \oplus H^{1}\left(  D^{n},\mathcal{R}%
^{0}\varphi\left(  \mathcal{T}^{0}\right)  \right) \label{ct}%
\end{equation}
The sheaves $\mathcal{R}^{q}\varphi\left(  \mathcal{T}^{0}\right)  ,q=0,1$ are
coherent by Grauert's Theorem, hence the terms $H^{1}$ and $H^{2}$ vanish in
the right-hand side. Therefore we have
\[
H^{1}\left(  X,\mathcal{T}^{0}\right)  \cong H^{0}\left(  D^{n},\mathcal{R}%
^{0}\varphi\left(  \mathcal{T}^{0}\right)  \right)  .
\]
The sheaf $\mathcal{R}^{1}\varphi\left(  \mathcal{T}^{0}\right)  $ vanishes at
the boundary of the polydisk since of (i). This yields that $H^{1}\left(
X,\mathcal{T}^{0}\right)  $ has finite dimension. Finally $\tau\doteq
\dim_{\mathbb{C}}T^{1}(X)<\infty$.

The second condition of Theorem \ref{2.1} can be weaken as follows

(ii)'$\,\mathcal{R}^{1}\varphi\left(  \mathcal{T}^{0}\right)  |\partial D^{n}=0.$

\begin{theorem}
\label{2.3} Let $(X,\varphi)$ be a polyhedron as in the previous theorem and
$f:\mathcal{X}\rightarrow(S,\circ)$ be its miniversal deformation. Then there
exists a neighborhood $S^{\prime}$ of the marked point $\circ$ in $S$ and a
closed subspace $M\subset S^{\prime}$ such that: \newline \textrm{(i)} for any
$s\in M$ the germ $(M,s)$ is maximal modular for the deformation $f$ of the
fiber $X_{s}$. \newline \textrm{(ii)} the restriction mapping $T^{0}%
(\mathcal{X}\times_{S}M/M)\rightarrow T^{0}(X)$ is surjective and $M$ is
maximal with this property; \newline \textrm{(iii)} the Zariski tangent space
$T_{\circ}(M)$ is the space of all vectors $t\in T_{\circ}\left(  S\right)  $
satisfying the equation $[\mathrm{D}_{\circ}f\left(  t\right)  ,v]=0$ for any
$v\in T^{0}(X)$; \newline \textrm{(iv)} the support of $M$ is the set of
points $s\in S^{\prime}$, where the mapping $\mathrm{D}_{s}f:T_{s}\left(
S\right)  \rightarrow T^{1}\left(  X_{s}\right)  $ is injective.
\end{theorem}

The statement (i) is essential; it means that any point of $M^{\prime}$ can be
taken as marked one, hence it is a globally defined complex analytic space. A
proof can be done on the lines of \cite{P3} and \cite{P6}. The parts
(ii)\label{1} and (iv) then follow from Theorem \ref{t1}.(i). The part (iii)
is a corollary of Theorem \ref{t1}.(iii), since the bracket $[\mathrm{D}%
f\left(  t\right)  ,v]$ is for any $t\in T_{\circ}\left(  S\right)  $ the
first obstruction for extension of $v$ to a tangent field on $\mathcal{X}%
\times_{S}M/M.$ The statement (iv) implies that \newline (\textrm{v})
\textit{if} $S$ \textit{is regular, then }$\mathrm{supp\,}M$ \textit{is given
by the equation} $\dim T^{1}\left(  X_{s}\right)  =\dim S,\,s\in S^{\prime}%
$.\newline \textbf{Definition 6.} Take a versal deformation $\left(
\mathcal{X},S\right)  $ as in Theorem \ref{2.1} and the maximal modular
stratum $M\subset S$ as in Theorem \ref{2.3}. We call $f|M:\mathcal{X}%
\times_{S}M\rightarrow M$ maximal modular deformation. Let $\mathcal{M}$ be a
complex analytic space and $Y$ be an analytic polyhedron over $\mathcal{M}$;
we call a \textit{modular }family the mapping $\mathcal{Y}\rightarrow
\mathcal{M}$, if it is maximal modular deformation of each fiber $Y_{s}%
,\;s\in\mathcal{M}$.

\section{Shrinking a polyhedron}

Let $\left(  g,q\right)  :\left(  Y,\psi\right)  \rightarrow\left(
X,\varphi\right)  $ be a morphism of analytic polyhedra. We call it
\textit{imbedding}, if $g$ is an imbedding, $q\ $is a coordinate projection
and there exists a proper holomorphic mapping $\rho:X\rightarrow\mathbb{C}%
^{n}$ that makes the commutative diagram%
\[%
\begin{tabular}
[c]{lll}%
$Y$ & $\overset{g}{\rightarrow}$ & $X$\\
$\downarrow\psi$ & $\swarrow\rho$ & $\downarrow\varphi$\\
$\mathbb{C}^{n}$ & $\overset{q}{\rightarrow}$ & $\mathbb{C}^{m}$%
\end{tabular}
\]
In geometric terms, this means that $\varphi_{i}=q_{i}\rho$ for $i=1,...,m$
and
\[
Y=\left\{  x\in X,\left|  \psi_{m+1}\left(  x\right)  \right|  <1,...,\left|
\psi_{n}\left(  x\right)  \right|  <1\right\}
\]
for an appropriate numeration of coordinates in $\mathbb{C}^{n}.$

\begin{proposition}
\label{tt}Let $\left(  g,q,\rho\right)  $ be an imbedding of polyhedra such
that \newline (i) $\mathrm{supp}\,\mathcal{T}^{1}\left(  X\right)  \Subset
g\left(  Y\right)  $ and \newline (ii) the mapping
\[
\rho:X\backslash g\left(  Y\right)  \rightarrow D^{m}\times\mathbb{C}%
^{n-m}\setminus D^{n-m}%
\]
is finite. Then we have $T^{1}\left(  Y\right)  \cong T^{1}\left(  X\right)  .$
\end{proposition}

\textsc{Proof}. Compare cohomology of tangent sheaves in $X$ and $Y.$ Consider
the diagram
\[%
\begin{tabular}
[c]{lllll}%
$H^{1}\left(  Y,\mathcal{T}^{0}\left(  Y\right)  \right)  $ &
$\cong$ & $H^{0}\left(  D^{n},\mathrm{R}^{1}\psi\left(
\mathcal{T}^{0}\left( Y\right)  \right)  \right)  $ &
$\overset{i}{\rightarrow}$ & $H^{0}\left(
D^{m},\mathrm{R}^{0}q\mathrm{R}^{1}\psi\left(
\mathcal{T}^{0}\left(
Y\right)  \right)  \right)  $\\
$\downarrow h_{1}$ &  &  &  & $\downarrow\gamma$\\
$H^{1}\left(  X,\mathcal{T}^{0}\left(  X\right)  \right)  $ &
$\cong$ & $H^{0}\left(  D^{m},\mathrm{R}^{1}\varphi\left(
\mathcal{T}^{0}\left( X\right)  \right)  \right)  $ &
$\overset{j}{\rightarrow}$ & $H^{0}\left(
D^{m},\mathrm{R}^{0}q\mathrm{R}^{1}\rho\left(
\mathcal{T}^{0}\left(
X\right)  \right)  \right)  $%
\end{tabular}
\]
The left-side isomorphisms follow from (\ref{ct}). The Leray spectral sequence%
\[
H^{r}\left(  D^{m},\mathrm{R}^{s}q_{\ast}\left(
\mathcal{F}\right)  \right)
\Rightarrow H^{r+s}\left(  D^{n},\mathcal{F}\right)  ,\;\mathcal{F}%
\doteq\mathrm{R}^{1}\psi\left(  \mathcal{T}^{0}\left(  Y\right)
\right)
\]
yields the isomorphism $i,$ since the sheaf $\mathcal{F}$ is
coherent. We have $\mathrm{R}^{r}q\,\mathrm{R}^{0}\rho\left(
\mathcal{T}^{0}\left(  X\right) \right)  =0$ for $r=1,2$ since the
sheaf $\mathrm{R}^{0}\rho\left(
\mathcal{T}^{0}\left(  X\right)  \right)  $ is coherent in $D^{m}%
\times\mathbb{C}^{n-m}$ and $q$ is a coordinate projection. By
applying the Leray sequence to the composition $\varphi=q\rho$ we
obtain the isomorphism
in $D^{n}$%
\[
\mathrm{R}^{1}\varphi\left(  \mathcal{T}^{0}\left(  X\right)
\right) =\mathrm{R}^{0}q\mathrm{R}^{1}\rho\left(
\mathcal{T}^{0}\left(  X\right) \right)
\]
This yields the bijection $j.$ Further we have
\[
\mathrm{R}^{1}\psi\left(  \mathcal{T}^{0}\left(  Y\right)  \right)
=\mathrm{R}^{1}\rho g\left(  \mathcal{T}^{0}\left(  Y\right)
\right) =\mathrm{R}^{1}\rho\,\mathrm{R}^{0}g\left(
\mathcal{T}^{0}\left(  Y\right) \right) =\mathrm{R}^{1}\rho\left(
\mathcal{T}^{0}\left(  X\right) |g\left(  Y\right)  \right)
\]
since $\mathrm{R}^{1}g=0.$ On the other hand, the sheaf $\mathrm{R}^{1}%
\rho\left(  \mathcal{T}^{0}\left(  X\right)  \right)  $ vanishes
in $X\backslash g\left(  Y\right)  $ since $\rho$ is finite there.
Therefore the right-hand side is isomorphic to the sheaf
$\mathrm{R}^{1}\rho\left( \mathcal{T}^{0}\left(  X\right)  \right)
.$ This yields the bijection
$\gamma$ and also the bijection $h_{1}.$ Consider the commutative diagram%
\[%
\begin{tabular}
[c]{lllll}%
$\mathrm{Ker}\,d_{2}\left(  Y\right)  $ & $\rightarrow$ &
$H^{0}\left( Y,\mathcal{T}^{1}\left(  Y\right)  \right)  $ &
$\overset{d_{2}\left( Y\right)  }{\rightarrow}$ & $H^{2}\left(
Y,\mathcal{T}^{0}\left(  Y\right)
\right)  $\\
$\downarrow k_{2}$ &  & $\downarrow h_{0}$ &  & $\downarrow h_{2}$\\
$\mathrm{Ker}\,d_{2}\left(  X\right)  $ & $\rightarrow$ &
$H^{0}\left( X,\mathcal{T}^{1}\left(  X\right)  \right)  $ &
$\overset{d_{2}\left( X\right)  }{\rightarrow}$ & $H^{2}\left(
X,\mathcal{T}^{0}\left(  X\right)
\right)  $%
\end{tabular}
\]
The isomorphism $h_{2}$ is defined similarly to $h_{1}$ and the bijective
mapping $h_{0}$ follows from (i). This diagram yields the isomorphism $k_{2}.$
This together with $h_{1}$ gives%
\[
T^{1}\left(  Y\right)  \cong H^{1}\left(  Y,\mathcal{T}^{0}\left(
Y\right) \right)  \oplus\mathrm{Ker}\,d_{2}\left(  Y\right)  \cong
H^{1}\left( X,\mathcal{T}^{0}\left(  X\right)  \right)
\oplus\mathrm{Ker}\,d_{2}\left( X\right)  \cong T^{1}\left(
X\right) .\blacktriangleright
\]

\begin{theorem}
\label{6.2} Let $\left(  Y,\psi\right)  ,\left(  X,\varphi\right)  $ be
analytic polyhedra that fulfil the conditions of Theorem \ref{2.1} and
$\left(  \mathcal{X},\Phi,S_{X}\right)  $ and $\left(  \mathcal{Y,}\Psi
,S_{Y}\right)  $ be their versal deformations. Let $\left(  g,q,\rho\right)
:\left(  Y,\psi\right)  \rightarrow\left(  X,\varphi\right)  $ be a imbedding
of analytic polyhedra that fulfil the conditions (i), (ii) of Proposition
\ref{tt}. Then there exists an imbedding of relative polyhedra $\left(
G,Q,R\right)  :\left(  \mathcal{Y,}\Psi,S_{Y}^{\prime}\right)  \rightarrow
\left(  \mathcal{X},\Phi,S_{X}^{\prime}\right)  $ such that $G\times
\circ=g,\,Q\times\circ=q,$ $R\times\circ=\rho$ and the diagram commutes
\begin{equation}%
\begin{tabular}
[c]{lllllll}%
$Y$ & $\rightarrow$ & $\mathcal{Y}$ & $\overset{G}{\rightarrow}$ &
$\mathcal{X}$ & $\leftarrow$ & $X$\\
$\downarrow\psi$ &  & $\downarrow\Psi$ & $\swarrow R$ & $\downarrow\Phi$ &  &
$\downarrow\varphi$\\
$\mathbb{C}^{n}\times\circ$ & $\subset$ & $\mathbb{C}^{n}\times S_{Y}^{\prime
}$ & $\overset{Q}{\rightarrow}$ & $\mathbb{C}^{m}\times S_{X}^{\prime}$ &
$\supset$ & $\mathbb{C}^{m}\times\circ$\\
& $\searrow$ & $\downarrow$ &  & $\downarrow$ & $\swarrow$ & \\
&  & $\mathbb{C}^{n}$ & $\overset{q}{\rightarrow}$ & $\mathbb{C}^{m}$ &  &
\end{tabular}
\label{di}%
\end{equation}
where $S_{X}^{\prime}$ and $S_{Y}^{\prime}$ are some neighborhoods of marked
points in $S_{X},$ respectively in $S_{Y}.$
\end{theorem}

\textsc{Sketch of the proof.} We can write $\rho=\varphi\times\xi$ where
$\xi:X\rightarrow\mathbb{C}^{n-m}$ is a bounded holomorphic mapping, say
$\left|  \xi\right|  <b$ in $\bar{X}$ for some positive $b.$ Consider the
polyhedron $\left(  X,\tilde{\rho}\right)  $ where $\tilde{\rho}=\varphi\times
b^{-1}\xi.$ It has the same boundary as $X$ and fulfils the conditions of
Theorem \ref{2.1}. Take a versal deformation $\left(  \mathcal{X},R,S\right)
$ of this polyhedron, where $R:\mathcal{X}\rightarrow\mathbb{C}^{n}\times S$
is the barrier function. We can write $R=R_{m}\times R_{n-m}\times F,$ where
$R_{m}:\mathcal{X}\rightarrow\mathbb{C}^{m},R_{n-m}:\mathcal{X}\rightarrow
\mathbb{C}^{n-m},F:\mathcal{X}\rightarrow S,$ and define the mapping
$\Phi=R_{m}\times F:\mathcal{X}\rightarrow\mathbb{C}^{m}\times S.$ We have
$R_{m}\times\circ=\varphi$ and the triple $\left(  \mathcal{X},R_{m},S\right)
$ is a deformation of the polyhedron $\left(  X,\varphi\right)  .$ Set
$\Psi=R_{m}\times bR_{n-m}\times F:\mathcal{X}\rightarrow\mathbb{C}^{n}\times
S;$ the polyhedron $\mathcal{Y}\doteq\left\{  x\in X\mathcal{;}\left|
bR_{n-m}\right|  <1\right\}  $ is embedded in $\mathcal{X}$ and $\mathcal{Y}%
\times\circ=Y,\Psi\times\circ=\psi,$ that $\left(  \mathcal{Y,}\Psi,S\right)
$ is a deformation of $\left(  Y,\psi\right)  .$ The diagram (\ref{di})
commutes, where $S_{X}^{\prime}=S_{Y}^{\prime}=S$ and $Q=q\times
\mathrm{id}_{S}.$ Analyzing the construction of $\left(  \mathcal{X}%
,R,S\right)  $ (see \cite{P6}) we can see that $\left(  \mathcal{X}%
,R_{m},S\right)  $ is isomorphic to a versal deformation $\left(
\mathcal{X},\Phi,S_{X}\right)  $ of the polyhedron $\left(  X,\varphi\right)
$. Next we prove that $\left(  \mathcal{Y,}\Psi,S\right)  $ is a versal
deformation of $\left(  Y,\psi\right)  $ and Theorem follows.
$\blacktriangleright$\newline \textbf{Definition. }The class of imbeddings in
$\left(  X,\varphi\right)  $ that satisfy the conditions of Proposition
\ref{tt} is called the \textit{germ} of $\left(  X,\varphi\right)  ;$ each
polyhedron $\left(  Y,\psi\right)  $ is a representative of this germ. By the
above theorem the versal deformations of all representative are naturally
isomorphic. So are the maximal modular strata $M$ and maximal modular deformations.

\section{Amalgams and automorphisms}

\noindent\textbf{Amalgam of modular deformations. }Let $f:\mathcal{X}%
\rightarrow M$ and $g:\mathcal{Y}\rightarrow N$ be modular families. If there
are points $s\in M,t\in N$ such that the fibers $X_{s}$ and $Y_{t}$ are
isomorphic, then there exist neighborhoods $U$ of $s$ and $V$ of $t$ and a
uniquely defined commutative diagram%
\[%
\begin{tabular}
[c]{lll}%
$f^{-1}\left(  U\right)  $ & $\overset{A}{\cong}$ & $g^{-1}\left(  V\right)
$\\
$\downarrow$ &  & $\downarrow$\\
$U$ & $\overset{a}{\cong}$ & $V$%
\end{tabular}
\]
This statement follows from Theorem \ref{6.2}. The above families can be
patched together along $a$ giving rise to the modular family $\mathcal{X}%
\sqcup_{a}\mathcal{Y}$ with the base $M\sqcup_{A}N,$ which is the amalgam
(coproduct) of $M$ and $N$:
\[%
\begin{tabular}
[c]{lllll}%
$f^{-1}\left(  U\right)  $ &  & $\;\overset{A}{\Longleftrightarrow}$ &  &
$g^{-1}\left(  V\right)  $\\
$\downarrow$ &  &  &  & $\downarrow$\\
$\mathcal{X}$ & $\rightarrow$ & $\mathcal{X}\sqcup_{a}\mathcal{Y}$ &
$\leftarrow$ & $\mathcal{Y}$\\
$\downarrow$ &  & $\;\downarrow$ &  & $\downarrow$\\
$M$ & $\overset{\mu}{\rightarrow}$ & $M\sqcup_{A}N$ & $\overset{\nu
}{\leftarrow}$ & $N$%
\end{tabular}
.
\]
where $\mu$ and $\nu$ are the natural morphisms. For given $f$ and $g$ several
local isomorphisms $\left(  a,A\right)  $ may occur which can make the amalgam
a complicated occasionally non-Hausdorff space. This is just the case, if the
automorphism group of $f$ or of $g$ is non-trivial.\newline
\textbf{Automorphisms.} Let $\left(  f,i\right)  $ be a miniversal deformation
of an analytic polyhedron $X$ as above and $M$ be the maximal modular stratum
in its base $S$. Let $X_{K}$ be the germ of $X$ on the compact set
$K\doteq\operatorname*{supp}\mathcal{T}^{1}\left(  X\right)  \Subset X$. We
show that an arbitrary automorphism $a$ of $X_{K}$ generates an automorphism
of the stratum $M$. Choose a polyhedron $X^{\prime}\subset X$ that contains
the germ $X_{K}$ such that $a$ defines a morphism $a:X^{\prime}\rightarrow X$.
Take a miniversal deformation $\left(  f^{\prime},i^{\prime}\right)  $ of
$X^{\prime}$; let $S^{\prime}$ be the base of $f^{\prime}$. By Theorem
\ref{6.2} there exists an imbedding $j:S^{\prime}\rightarrow S$ that induces
an isomorphism $f\times j\cong f^{\prime}$ and maps the modular stratum
$M^{\prime}$ to the modular stratum $M$. Consider the deformation $(f,ia)$.
This is a miniversal deformation of $X^{\prime},$ hence it is induced from
$(f^{\prime},i^{\prime})$ by a germ endomorphism $\alpha:S\rightarrow
S^{\prime}$. It is an automorphism, since $\mathrm{D}\alpha$ is bijective at
the marked point. The restriction $\left.  \alpha\right|  M$ is uniquely
defined and the composition $\alpha_{M}\doteq j^{-1}\left.  \alpha\right|
M:M\rightarrow M$ is an automorphism of the germ $M$.

\begin{corollary}
\label{c2}The mapping $a\mapsto\alpha_{M}$ defines a homomorphism of the group
$\mathrm{Aut}\left(  X_{K}\right)  $ to the group $\mathrm{Aut}\left(
M\right)  $ of automorphisms of the maximal modular germ $M.$
\end{corollary}

The kernel $\mathrm{Aut}_{0}\left(  X_{K}\right)  $ of this homomorphism is a
normal subgroup and the quotient $\mathrm{G}\left(  X_{K}\right)
=\mathrm{Aut}\left(  X_{K}\right)  /\mathrm{Aut}_{0}\left(  X_{K}\right)  $
acts faithfully; we call it \textit{active} automorphism group. Any
automorphism $a$ of $X_{K}$ generated by a tangent field $v\in T^{0}\left(
X_{K}\right)  $ acts trivially on $M.$ The field $v$ generates a local
holomorphic subgroup $a\left(  \zeta\right)  $ of automorphisms. The
infinitesimal action of this group in $T_{\circ}(S)$ is given by the
commutator $t\mapsto\left[  \mathrm{D}_{\circ}f\left(  t\right)  ,v\right]  $,
which is trivial in $M$ due to Theorem \ref{2.3}.(iii). Therefore
$a_{M}=\mathrm{id}$ that is $a\in$ $\mathrm{Aut}_{0}\left(  X_{K}\right)  $.
The active group \textrm{G}$\left(  X_{K}\right)  $\textrm{\ }is discontinuous
in several examples (see below).

The fibers $X_{s}$ are, of course, isomorphic for points $s$ in any coset of
the automorphism group. The inverse statement is by no means obvious.\newline
\textbf{Problem 1} Let $f:\mathcal{X}\rightarrow\left(  M,\circ\right)  $ be a
maximal modular deformation of a space $X\cong X_{\circ}$. When does the
existence of a isomorphism $X_{s}\cong X_{t}$ for some points $s,t\in M$ imply
$\alpha_{M}\left(  s\right)  =t$ for an element $a\in\mathrm{G}\left(
X_{s}\right)  $? In particular, when an isomorphism $X_{s}\cong X,\,s\in
M^{\prime}$ implies $s=\circ$?

Is it true, if the group $\mathrm{G}(X_{s})$ is finite?

Let $\mathcal{X}\rightarrow M$ be a modular family with a irreducible base $M$
and $s\in M$. By Corollary \ref{c2} any element $\alpha\in\mathrm{G}\left(
X_{s}\right)  $ of the fiber $X_{s}=f^{-1}\left(  s\right)  $ generates an
automorphism $\alpha_{M}$ of the germ $\left(  M,s\right)  $. By analytic
continuation, we can extend $\alpha_{M}$ to a uniquely defined automorphism
$\alpha$ of the base $M$. Let $\mathrm{G}$ be automorphism group of $M$
generated by all the elements $\alpha_{M}$ for $\alpha\in\mathrm{G}\left(
X_{s}\right)  ,s\in M$. The quotient $M/\mathrm{G}$ can be taken as a
candidate for `true' moduli space.

\section{Modular deformation of compact spaces}

We start with very classical families that appear to be modular.\newline
\textbf{Example 1. }The family $T=\left\{  T\left(  \lambda\right)  \right\}
,\,\lambda\in\mathbb{C}_{+}$ of $1$-tori parameterized by the upper half-plane
$\mathbb{C}_{+}$; $T\left(  \lambda\right)  $ is the quotient of the plane
$\mathbb{C}$ by the lattice generated by $1$ and $\lambda.$ This family is
maximal modular for each point $\lambda\in\mathbb{C}_{+}$. The automorphism
group of the family is the group $\mathrm{Sl}\left(  2,\mathbb{Z}\right)  $
which acts in $\mathbb{C}_{+}$ by fractional linear transformations. The group
$\mathrm{Aut}_{0}\left(  T\left(  \lambda\right)  \right)  $ is equal to the
semi-direct product of the group $T\left(  \lambda\right)  $ generated by
tangent fields on $T\left(  \lambda\right)  \mathbb{.}$ The group
$\mathbb{Z}_{2}$ generated by the mapping $z\mapsto-z$ in $\mathbb{C}$ acts
trivially in the base $\mathbb{C}_{+}$. Take the standard fundamental domain
$D=\left\{  \lambda:\left|  \mathrm{Re\,}\lambda\right|  <1/2,\left|
\lambda\right|  >1\right\}  $ and consider the restriction of the family to
$D.$ There are two points $\lambda_{2}=\imath,\lambda_{3}=\sqrt[3]{-1}$ on the
boundary of $D$ such that the active automorphism group $\mathrm{G}(T\left(
\lambda\right)  )$ has elements with non-trivial action: these are elements
$a_{j}\in\mathrm{G}\left(  T\left(  \lambda_{j}\right)  \right)  $ of order
$2$, respectively $3.$ These elements are generated by rotation of the unit
square by $\pi/2$ and by rotation of the rhombus by $\pi/3,$ respectively.
According to previous section, these groups generate transformation group
$\mathrm{G}$ of $D.$ The quotient space $D/\mathrm{G}$ is isomorphic to the
complex plane $\mathbb{C}$ and the family $T$ gives rise to a family
$\tilde{T}$ of tori on $\mathbb{C}.$ It can be compacted by means of amalgam
with the deformation $f:\mathcal{Y}\rightarrow U$ of the singular curve
$Y_{0},$ which is the projective line with one point of transversal
self-intersection. Here $U$ is the unit disk and the deformation is given in
an affine chart by
\[
w^{2}-z(z-s)(z-1)=0,\;\left|  s\right|  <1/2
\]
The fibers $Y_{s}$ are non-isomorphic tori for $s\neq0$, and $\dim
T^{1}\left(  Y_{0}\right)  =1,$ hence the family $\mathcal{Y}$ is maximal
modular. It can be patched to $\tilde{T}\rightarrow\mathbb{C}$ giving rise to
a modular family with base $\mathbb{C}\cup\{s=0\}=\mathbb{CP}^{1}$. The point
$s=0$ corresponds to infinity, since a ratio of periods of the surface $Y_{s}$
tends to infinity as $s\rightarrow0$. \newline \textbf{Example 2.} Generalize
the above construction for curves of an arbitrary genus $g>1$. Consider the
family of hyperelliptic Riemann surfaces $X_{a},\;a\in\mathbb{C}^{m}$, where
$X_{a}$ is given in the affine chart by
\begin{equation}
w^{2}-p(z)=0,\;p(z)=z^{m}+a_{1}z^{m-1}+...+a_{m},\;m=2g+2. \label{hyp}%
\end{equation}
Suppose that $g\geq4$, take the singular surface $X_{0}$ for which
$a_{1}=...=a_{m}=0$ and check that%

\begin{equation}
\dim T^{1}(X_{0})=3g-3 \label{t}%
\end{equation}
We have
\[
\dim T^{1}(X)=\dim H^{0}(X,\mathcal{T}^{1})+\dim H^{1}(X,\mathcal{T}^{0}),
\]
since $H^{2}(X,\mathcal{T}^{1})=0$. Further,
\[
\dim H^{0}(X,\mathcal{T}^{1})=\dim\mathcal{T}_{s}^{1}=m-1=2g+1,
\]
since the only singular point $s=(0,0)$ is of multiplicity $m-1$. The surface
$X$ is union of two spheres that are tangent one to another at the origin. The
tangent sheaf $\mathcal{T}^{0}$ is generated at the origin by two fields
\[
t_{1}=mz\frac{\partial}{\partial z}+2w\frac{\partial}{\partial w}%
,\;t_{2}=2w\frac{\partial}{\partial z}-mz^{m-1}\frac{\partial}{\partial w}%
\]
over the algebra over $\mathcal{O}(X)_{s}$. It helps to check that $\dim
H^{1}(X,\mathcal{T}^{0})=g-4$, which implies (\ref{t}).

Let $F_{0}:\mathcal{Y}\rightarrow S$ be a miniversal deformation of $X_{0}$.
Prove that it is universal. The base $S$ is a piece of $\mathbb{C}^{3g-3}$ and
any non-singular fiber $Y_{s}$ is a Riemann surface of genus $g$. Therefore,
we have again $\dim T^{1}(Y_{s})=\dim H^{1}(Y_{s},\mathcal{T}^{0})=3g-3$. This
implies that the dimension of $T^{1}(Y_{s})$ is constant in the deformation
$F,$ since it is a upper semi-continuous function on the base. By Theorem
\ref{2.3}.(iv) the deformation $F$ is maximal modular and therefore universal.
The family of surfaces $X_{a}$ is obviously a deformation of $X_{0}$ which
yields that $\mathrm{dim}\,T^{1}(X_{a})=3g-3$. Therefore any miniversal
deformation $F_{a}$ of $X_{a}$ is universal too. They can be amalgamated
together in a maximal modular family $F$ with a base $S$ which is an open
subset in $\mathbb{C}^{3g-3}$. On the other hand, $S$ contains the base
$\mathbb{C}^{m}$ of the family (\ref{hyp}). The set $D\subset S$ of critical
values of this deformation coincides with the variety $\Delta(p)=0$ in
$\mathbb{C}^{m}$, where $\Delta(p)$ is the discriminant of the polynomial $p$.

On the other hand, the Teich\"{m}uller space $T(g)$ is the base of a universal
family $R(g)$, whose fibers are all non-singular Riemann surfaces of genus $g$
with an additional structure. The automorphism group $\Gamma(g)$ of the family
$R(g)$, called modular group, acts discontinuously in $T(g)$. The deformation
$F$ is amalgamated with the family $R(g)$ and the base $S\backslash D$ is
patched to $T(g)/\Gamma(g)$. This gives a compacting of the space
$T(g)/\Gamma(g)$, whereas the discriminant set $D$ covers the boundary.

In the cases $g=2,3$ the deformation $F$ as above is no more modular for the
surface $X_{0}$, since $\mathrm{\dim\,}T^{1}(X_{0})>3g-3$. We take the surface
$X(q)$ for $q(z)=(z^{2}-1)z^{2g}$ instead. The miniversal deformation of
$X(q)$ is universal, like in the case $g=1$.

This compacting of the space $T(g)$ is, apparently, different from that of
W.Baily \cite{B}.\newline \textbf{Example 3.} Fix an integer $n>1$ and
consider a proper holomorphic mapping of manifolds $F:T\rightarrow S$ whose
fibers are complex analytic $n$-tori, see \cite{KS}. The dimension of $S$ is
equal to $n^{2}$ and the mapping is maximal modular deformation of each fiber
as in the case $n=1.$ On the other hand, the active group $\mathrm{G}$ of $F$
is not discontinuous. Moreover, any non-empty open set $U$ contains a point
$s$ such that the coset $\mathrm{G}s\cap U$ is infinite. The moduli space
$M/\mathrm{G}$ is not separable.\newline \textbf{Example 4.} Take the family
$F:\mathcal{V}\rightarrow P(m)$ of algebraic hypersurfaces of degree $m$ in
$\mathbb{CP}^{n},\;n\geq3$, where $P(m)$ is the space of all homogeneous
polynomials of degree $m$ and $V(f)\doteq F^{-1}(f)$ is the hypersurface
defined by the equation $f=0$ in $\mathbb{CP}^{n}$. Take a polynomial $f\in
P(m)$ such that the variety $V(f)$ is non-singular and choose an affine
subspace $S\subset P(m)$ that contains the point $f$ and is transversal to the
linear span of polynomials $z_{i}\partial f/\partial z_{j}$, $i,j=0,...,n.$
Consider the restriction $\left.  F\right|  S$ of this family; the
Kodaira-Spencer mapping $\left.  DF\right|  S:T_{f}(S)\rightarrow T^{1}(V(f))$
is injective. This mapping is surjective and the family $\left.  F\right|  S$
is versal, except in the case $n=3,\,m=4$, see \cite{KS}. Moreover, this
family is maximal modular, since $\dim S=(m+n)!/n!\,m!-(n+1)^{2}$ does not
depend on $f$. The automorphism group of any fiber $X$ is finite for the same cases.

In the exceptional case $n=3,\,m=4$ the miniversal deformation $F$ of any
surface contains non-algebraic fibers $X_{s}$ that are also $\mathrm{K3}%
$-surfaces (the canonical bundle is trivial). The base $S$ has dimension $20$
and the algebraic fibers only appear for points $s$ in the dense subset
$S_{alg}$ which is a countable union of $19$-dimensional subspaces in the
base, see more details in \cite{KS}. Therefore $\dim T^{1}(X_{s})=20$ on the
germ $S$ and the formation $F$ is again maximal modular. On the other hand,
there is no non-trivial tangent fields on a $\mathrm{K3}$-surface, but its
automorphism group can be infinite \cite{MM}.

\section{Modular deformations of singular points}

The modular deformations of polyhedra with isolated singularities look similar
to the above examples. Moreover, there is a relation between deformations in
these two categories.

Let $D^{n}$ be the open unit polydisk in a coordinate space $\mathbb{C}^{n}$
centered at the origin. Take holomorphic functions $f_{1},...,f_{k}$ in
$\bar{D}^{n}$ and consider the analytic polyhedron $X=\{z\in D^{n}%
,\,f_{1}(z)=...=f_{k}(z)=0\}$ endowed with the sheaf $\mathcal{O}%
(X)=\mathcal{O}(D_{\varepsilon}^{n})/\mathcal{I}$, where $\mathcal{I}$ is the
sheaf-ideal generated by these functions. The polyhedron $X$ satisfies the
condition \ref{2.2}, if it contains only finite number of singular points.

Suppose that $k=1$, the generating function $f$ is weighted homogeneous for
certain coordinate system $z_{1},...,z_{n}$ in $\mathbb{C}^{n}$ and there is
only one singular point $z=0$ in $X$ (otherwise the dimension of the singular
set of $X$ is positive). The base $M$ of the maximal modular deformation is
isomorphic to the subspace of $T^{1}\left(  X\right)  $ of elements whose
weight is equal to that of $f$. In particular, $M$ is a simple point for
singularities of types $A,D$ and $E$ that are characterized by the inequality
$weight\left(  \tau\right)  <weight\left(  f\right)  $ for all $\tau\in
T^{1}\left(  X\right)  $.

Assume now that the weights of the coordinates are equal to $1$ so that $f$ is
a homogeneous polynomial of degree $m>0$. It defines a non-singular
hypersurface $V(f)$ in $\mathbb{CP}^{n-1}$ and any deformation of $V(f)$ as in
Example 4 generates a deformation of the analytic polyhedron $X(f)\doteq\{z\in
D^{n},\,f(z)=0\}$. We have
\begin{equation}
T^{1}(X(f))\cong\mathbb{C}\left[  z_{1},...,z_{n}\right]  /(\partial
f/\partial z_{1},...,\partial f/\partial z_{n}), \label{iso}%
\end{equation}
The weight grading generates a grading in the space (\ref{iso}).

\begin{proposition}
Suppose that $n\geq4$ and the case $n=m=4$ is excluded. Then there exists a
linear injective mapping
\begin{equation}
j:T^{1}(V(f))\rightarrow T^{1}(X(f)), \label{emb}%
\end{equation}
whose image coincides with the subspace of weight $m$.
\end{proposition}

$\blacktriangleleft$ By diagram (9.3) of \cite{P2}, p.109 we have
$T^{1}(V(f))=P(m)/J$, where $P(m)$ is the space of polynomials in
$\mathbb{C}^{n}$ of degree $m$ and $J$ is the linear envelope of the
polynomials $z_{i}\partial f/\partial z_{j},\,i,j=1,...,n$. Note that the term
$H^{1}(V(f),\theta)$ vanishes in the diagram due to Proposition 9.2. The
isomorphism (\ref{iso}) makes the mapping (\ref{emb}) obvious.
$\blacktriangleright$

\begin{corollary}
For any family of homogeneous polynomials $\{f_{s},\,s\in S\}$ of $n$
variables with the only critical point $z=0$ the numbers \textrm{dim\thinspace
}$T^{1}\left(  V(f_{s})\right)  $ and $\mathrm{dim}\,T^{1}\left(
X(f_{s})\right)  $ stay constant. If one of the families is modular, another
is also modular.
\end{corollary}

$\blacktriangleleft$ For the first statement we note that $\mathrm{dim}%
\,T^{1}(X(f_{s}))=(m-1)^{n}.$ By the previous Proposition $\mathrm{\dim
\,}T^{1}(V(f_{s}))$ is equal to the number of monomials of degree $m$ minus
$n,$ that is the dimensions depend only on $n$ and $m.$ We have the equation
for the Kodaira-Spencer mappings $\mathrm{D}\Phi=j\mathrm{D}F$ where $F$ and
$\Phi$ are the families of polyhedrons and of projective hypersurfaces,
respectively. The second statement follows from injectivity of
$j.\;\blacktriangleright$

The modular stratum of a non weighted-homogeneous hypersurface may be singular
and have imbedded primary components.\newline \textbf{Example 5.} Consider the
polyhedron in $D^{3}$ defined by the equation
\begin{equation}
f_{p,q,r}\left(  \lambda;x,y,z\right)  \equiv x^{p}+y^{q}+z^{r}+\lambda xyz=0
\label{fl}%
\end{equation}
in the complex coordinates $x,y,z$. We denote the singular germ at the origin
by $T_{p,q,r}(\lambda)$, if the equation holds $1/p+1q+1/r=1$. The generating
function is then weighted homogeneous and the parameter $\lambda$ is the
coordinate in a maximal modular deformation, where $\lambda$ runs over the
complex plane with few gaps, where the polyhedron contains a singular curve.
This follows from Theorem \ref{2.3}.(iv), since the space of elements $t\in
T^{1}\left(  X_{\lambda}\right)  $ such that $\left[  t,v\right]  =0$ for all
$v\in T^{0}\left(  X_{\lambda}\right)  $ is one-dimensional. There are just
three modular families of this type: $T_{3,3,3}(\lambda),\,T_{4,4,2}(\lambda)$
and $T_{6,3,2}(\lambda)$ whose Tyurina numbers $\tau\left(  X\right)
\doteq\dim T^{1}\left(  X\right)  $ are $8,9,10,$ respectively.

In the case $1/p+1q+1/r<1$ the surfaces (\ref{fl}) are all isomorphic for
$\lambda\neq0$. We set $\lambda=1$ and use the notation $T_{p,q,r}$%
.\newline \textbf{Example 6.} Start with the family $T_{3,3,3}\left(
\lambda\right)  $ defined for $\lambda\in\mathbb{C}$. We have $\tau\left(
X\right)  =8$ and the family is maximal modular for $\lambda^{3}\neq-27$. The
active group $\mathrm{G}=\mathbb{Z}_{3}$ acts in the family by $\lambda
\mapsto\varepsilon\lambda$, where $\varepsilon^{3}=1$. Take the deformation
$Y\rightarrow N$ of the germ $T_{4,3,3}$ defined in $D^{3}\times N$ by the
function
\[
F(r,s,t;x,y,z)=x^{4}+y^{3}+z^{3}+\mu xyz+rx^{3}+s_{1}y+s_{2}y^{2}+t_{1}%
z+t_{2}z^{2},
\]
where the number $\mu\neq0$ is fixed. It is maximal modular for the base
$N\subset M\times{\mathbb{C}}^{3}$ given by the ideal $I=I_{1}\cap I_{0}$
where the ideal $I_{1}=(s_{1},s_{2},t_{1},t_{2})$ defines the line
parameterized by the coordinate $r$ and we have an isomorphism $Y_{r,0,0}%
\simeq T_{3,3,3}(\lambda)$ for $\lambda=r^{-1/3}$. The ideal $I_{0}$
determines a fat point at the origin in ${\mathbb{C}}^{5}$, whose embedding
dimension is equal to $5$. Therefore the modular deformations $T_{3,3,3}%
(\lambda)$ and $Y$ are patched together to a maximal modular deformation with
a compact base $M$.\newline \textbf{Example 7. } Consider the family
$\mathcal{X}$ of curves $X_{\lambda}\subset D_{1}^{2}$ given by the equation
\[
(x^{3}-xy^{2})(ax-by)=0,\;\lambda=b/a\in\mathbb{CP}^{1}.
\]
The family is maximal modular for all $\lambda$, except in three points
$\lambda=-1,0,1$. These gaps can be patched by means of another modular
deformations. B.Martin has considered the family $Y\rightarrow\mathbb{C}^{7}$
of polyhedra in $D^{2}\times\mathbb{C}^{7}$, given by the equation
\[
x^{4}-x^{2}y^{2}+s_{1}x+s_{2}y+s_{3}xy+s_{4}y^{2}+s_{5}y^{3}+s_{6}%
xy^{2}+ty^{4}+y^{5}=0.
\]
He has shown that this family is maximal modular over the base $M,$ whose
ideal is $I(M)\doteq J_{1}\cap J_{0}$, where $J_{1}=(s_{1},...,s_{6})$ defines
the line and the ideal $J_{0}$ has $9$ generators and defines a fat point at
the origin. For any point $(0,t)\in\mathbb{C}^{7}$, $t\neq0$ the germ
$Y_{0,t}$ is isomorphic to $X_{\lambda}$ for $\lambda=2\sqrt{t}$. Therefore
the family $\mathcal{Y}$ can be amalgamated with the family $\left\{
X_{\lambda}\right\}  $ patching the gap $\lambda=0$. The gaps $\lambda=\pm1$
are patched in a similar way. The resulting amalgam is a maximal modular
family $\tilde{T}_{3,3,3}$ with the compact base $\mathbb{CP}^{1}$.

Note that the Milnor number of a fiber of this family is not constant; it is
equal to $9$, except in the points $\lambda=-1,0,1,$ where it equals
$10.$\newline \textbf{Example 8. }Consider the family $\mathcal{Z}$ of
complete intersection curves $Z_{s}\subset D^{3}\times\mathbb{C}$ defined by
two polynomials:
\[
x^{4}+y^{4}+2z^{2}=sz-xy=0,\;s\in\mathbb{C.}%
\]
It is modular with $\tau\left(  Z_{s}\right)  =9$, \cite{Al1}. The curve
$Z_{s}$ is isomorphic for $s\neq0$ to the plane curve $x^{4}+y^{4}%
+2s^{-2}x^{2}y^{2}=0$. The germ at the origin of this curve is isomorphic to
the fiber $X_{\lambda}$ of the family in Example 7 with $\lambda=s^{2}$. Now
the family $\mathcal{Z}$ is amalgamated with $\mathcal{X}$ mending the gap at
the point $\lambda=0$. This provides a compacting of the family $\mathcal{X}$
at this point, which is different form that of Example 7.\newline
\textbf{Example 9. }(\cite{HM}) The maximal modular deformations of polyhedra
$Y_{7,3,2},\;Y_{6,4,2}$ and $Y_{6,3,3}$ can be amalgamated with the modular
family $T_{6,3,2}(\lambda)$. The modular deformation of $Y_{6,4,2}$ is given
by the function
\[
x^{6}+y^{4}+z^{2}+xyz+rx^{4}+sx^{5}+ty^{3}+vz=0
\]
in $D^{3}\times\mathbb{C}^{4}$. It is maximal modular over the base germ $M$
defined by the ideal $I(M)=J_{1}\cap J_{2}\cap J_{3},$ where
\[
J_{1}=(r,s,v),\;J_{2}=\left(  r,s^{2},t^{2},20v-st\right)  \;,J_{3}=\left(
t,4r-s^{2},v\right)  .
\]
The zero set of $J_{1}$ is the line with the coordinate $t$; the fiber
$Y(0,0,t,0)$, $t\neq0$ is a surface whose germ is isomorphic to $T_{6,3,2}%
(\lambda)$ for $\lambda=t^{-1/3}.$ The ideal $J_{2}$ defines the fat point at
the origin and $J_{3}$ defines a smooth curve parameterized by $s_{5}$; the
fiber of deformation over this curve is the union of the germ $T_{4,4,2}%
(\lambda)$ for some $\lambda$ and of a singular germ of multiplicity $1$ at
the point $(-s/2,0,0)$.

The maximal modular deformation $Y_{6,3,3}$ has the base, which is the union
of three components. Two of them are straight lines that are symmetric under
the permutation $x\leftrightarrow y$; the fibers are of type $T_{6,3,2}%
(\lambda)$. The third component is a curve; the general fiber of the
deformation over this curve has the singular germ of type $Y_{4,3,3}$ at the
origin and another singular point of multiplicity $2$.

The modular deformation $Y_{7,3,2}$ is similar to the previous two
cases.\newline \textbf{Problem 2.} These examples give rise to the general
questions: let $f^{\ast}$ be a modular deformation over the punctured disk
$D^{\ast}$. When there exists a modular deformation $f$ over $D$ such that
$f|D^{\ast}\cong f^{\ast}$? Note that the deformation $f$ need not to be
unique as we saw in Examples 8 and 9.

\section{Concluding remarks}

Let $X$ be an arbitrary hypersurface germ with isolated singularity which is
not weighted homogeneous. Then the tangent space $T\left(  M\right)  $ to the
maximal modular stratum has always positive embedding dimension. Indeed, the
action of the Lie algebra $T^{0}\left(  X\right)  $ in $T^{1}\left(  X\right)
$ is nilpotent and by Engel's theorem, there is an element $t\in T^{1}\left(
X\right)  $, $t\neq0$ that vanishes under action of the Lie algebra. By
Theorem \ref{2.3} this yields $t\in T\left(  M\right)  $. This is the case for
any singular germ of type $T_{p,q,r}$.

More examples can be extracted from the classification from Arnold's list
\cite{AGV}, where several families of singularities with constant Milnor
number are listed. Such a family has also constant Tyurina number, except for
a proper subvariety of the family base, where this number jumps up.
Restricting the family to a certain affine subvariety $M$, yields a modular
deformation. The `gaps' in $M$ could be filled by means of amalgams with
appropriate modular deformations. The completed variety $\tilde{M}$ is
expected to be a projective algebraic. Besides, any `solitary' singularity in
Arnold's list which is not included in families may have the maximal modular
deformation with the base $M$ which is either a fat point or a germ of
positive dimension with splitting singularity like that of $T_{6,4,2}$ and
$T_{6,3,3}$ in Example 9.

A.Aleksandrov \cite{Al} has found a series of modular families of complete
intersection curves in $\mathbb{C}^{3}$. Examples of families with constant
Tyurina number $\geq34$ were calculated in the book \cite{LP}. B.Martin
\cite{Ma2}. T.Hirsch and B.Martin \cite{HM} are found sophisticated examples
of modular deformations.

\end{document}